\documentclass[10pt,a4paper]{amsart}

\usepackage{amssymb}
\usepackage{amsthm}
\usepackage{amsfonts}
\usepackage{microtype}
\usepackage{mathrsfs}
\usepackage{graphicx}
\usepackage{color}
\usepackage{latexsym}
\usepackage{rotating}
\usepackage{xspace}
\usepackage[all]{xy}
\usepackage{longtable}
\usepackage{leftidx}
\usepackage{mathtools}
\usepackage{titlesec}
\usepackage{tikz}
\usetikzlibrary{calc}
\usetikzlibrary{arrows}
\usetikzlibrary{decorations.markings}
\usepackage{geometry}

\newtheorem{theorem}{Theorem}[section]
\numberwithin{equation}{theorem}
\newtheorem{lemma}[theorem]{Lemma}
\newtheorem{proposition}[theorem]{Proposition}
\newtheorem{corollary}[theorem]{Corollary}

\theoremstyle{definition}
\newtheorem{definition}[theorem]{Definition}
\newtheorem{example}[theorem]{Example}

\newtheorem{remark}[theorem]{Remark}

\theoremstyle{conjecture}

\newtheorem{question}[theorem]{Question}

\newcommand{\Ass}{\operatorname{Ass}}

\newcommand{\grade}{\operatorname{grade}}

\newcommand{\Assh}{\operatorname{Assh}}

\newcommand{\ara}{\operatorname{ara}}

\newcommand{\cd}{\operatorname{cd}}

\newcommand{\id}{\operatorname{id}}
\newcommand{\fd}{\operatorname{fd}}

\newcommand{\pd}{\operatorname{pd}}

\newcommand{\V}{\operatorname{V}}
\newcommand{\RH}{\operatorname{H}}

\newcommand{\Ext}{\operatorname{Ext}}
\newcommand{\Supp}{\operatorname{Supp}}

\newcommand{\Hom}{\operatorname{Hom}}

\newcommand{\Ann}{\operatorname{Ann}}
\newcommand{\Rad}{\operatorname{Rad}}

\newcommand{\Max}{\operatorname{Max}}

\newcommand{\lo}{\longrightarrow}
\newcommand{\fl}{\frak{l}}
\newcommand{\fm}{\frak{m}}

\newcommand{\fa}{\frak{a}}
\newcommand{\fb}{\frak{b}}

\newcommand{\fn}{\frak{n}}

\newcommand{\suchthat}{\;\ifnum\currentgrouptype=16 \middle\fi|\;}

\newenvironment{prf}[1][Proof]{\begin{proof}[\bf #1]}{\end{proof}}

\makeatletter
\newcommand{\holim@}[2]{%
\vtop{\m@th\ialign{##\cr
\hfil$#1\operator@font holim$\hfil\cr
\noalign{\nointerlineskip\kern1.5\ex@}#2\cr
\noalign{\nointerlineskip\kern-\ex@}\cr}}%
}
\newcommand{\holim}{%
\mathop{\mathpalette\holim@{\rightarrowfill@\textstyle}}\nmlimits@
}
\makeatother

\makeatletter
\def\@secnumfont{\bfseries}
\def\section{\@startsection{section}{1}%
\z@{.7\linespacing\@plus\linespacing}{.5\linespacing}%
{\normalfont\Large\bfseries\filcenter}}
\def\subsection{\@startsection{subsection}{2}%
\z@{.5\linespacing\@plus.7\linespacing}{-.5em}%
{\normalfont\large\bfseries}}
\makeatother
\begin{document}

\author[P. Pourghobadian, K. Divaani-Aazar and A. Rahimi]
{Parisa Pourghobadian, Kamran Divaani-Aazar and Ahad Rahimi}

\title[Relative Cohen-Macaulay modules ...]
{Relative Cohen-Macaulay modules under ring homomorphisms}

\address{P. Pourghobadian, Department of Mathematics, Faculty of Mathematical Sciences, Alzahra
University, Tehran, Iran.}
\email{paparpourghobadian@gmail.com}

\address{K. Divaani-Aazar, Department of Mathematics, Faculty of Mathematical Sciences, Alzahra
University, Tehran, Iran.}
\email{kdivaani@ipm.ir}

\address{A. Rahimi, Department of Mathematics, Razi University, Kermanshah, Iran.}
\email{ahad.rahimi@razi.ac.ir}



\begin{abstract} Let $R$ be a commutative Noetherian ring with identity (not necessarily local) and $\fa$ a proper ideal of
$R$. We study the invariance of some classes of $\fa$-relative Cohen-Macaulay modules under pure ring homomorphisms and ring
homomorphisms of finite flat dimension. Our results extend several results in the existing literature on homological modules.
\end{abstract}

\maketitle

\tableofcontents

\section{Introduction}

Throughout this article, the word ring stands for commutative Noetherian rings with identity. We shall refer to a ring (resp.
module) with a homological criterion as a homological ring (resp. homological module). The theory of homological rings was
started in 1954, when Auslander, Buchsbaum, and Serre proved a celebrated homological criterion for regular local rings.
The most important classes of homological rings satisfy the following implications:
$$\text{Regular\ ring} \Rightarrow  \text{complete\ intersection\ ring} \Rightarrow \text{Gorenstein\ ring}\Rightarrow
\text{Cohen-Macaulay\ ring}$$ $$ \Rightarrow \text{surjective Buchsbaum\ ring} \Rightarrow  \text{Buchsbaum\ ring}
\Rightarrow \text{quasi-Buchsbaum\ ring}$$ $$\Rightarrow \text{generalized\ Cohen-Macaulay\ ring}.$$ From the Gorenstein
rings on, the above classes of rings have module theoretic counterparts. These classes of modules make the core of the
theory of homological modules. Among them, Cohen-Macaulay modules have been one of the main streams of research in
commutative algebra for decades.

There are comprehensive results on the ascent and descent of homological properties along local ring homomorphisms. In
particular, it is known that if $f:R\lo T$ is a surjective local ring homomorphism such that $\ker f$ is generated
by an $R$-regular sequence, then the Cohen-Macaulay and Gorenstein properties ascend and descend from $R$ to $T$.
There have been a number of studies exploring the ascent and descent of homological properties along general local ring
homomorphisms of finite flat dimension; see e.g. \cite{AF} and \cite{AFH}. The purpose of this paper is to investigate
the ascent and descent of some classes of relative Cohen-Macaulay modules along ring homomorphisms.

There are several generalizations of the notion of Cohen-Macaulay modules in the literature, among them is the notion of
relative Cohen-Macaulay modules. This concept is introduced by Hellus and Schenzel \cite{HS}, and  Rahro Zargar and Zakeri
\cite{RZ2}. Let $\fa$ be a proper ideal of a ring $R$ and $M$ a finitely generated $R$-module. The $R$-module $M$ is called
$\fa$-{\it relative Cohen-Macaulay} if $\RH_{\fa}^{i}\left(M\right)=0$ for all $i\neq \grade(\fa,M)$. Several authors followed
up on relative Cohen-Macaulay modules; see e.g. \cite{HeSt, Sc1, Sc2, R, JR, Ra2, RZ1, CH, Ra1, DGTZ1, DGTZ2}. Particularly,
\cite{DGTZ2} introduces the notion of $\fa$-relative system of parameters, that seems to be quite useful for studying relative
Cohen-Macaulay modules.

In \cite{DGTZ1}, the authors present a relative theory of generalized Cohen-Macaulay, quasi-Buchsbaum, Buchsbaum, and
surjective Buchsbaum modules, by using the notion of $\fa$-relative system of parameters. Also, a relative theory of
maximal Cohen-Macaulay, Gorenstein and regular modules is presented in \cite{PDR}.

In this paper, we examine the ascent and descent of relative Cohen-Macaulayness, relative maximal Cohen-Macaulayness,
relative Gorensteinness, and relative regularness along ring homomorphisms.  The two types of ring homomorphisms that
we will examine are: Pure ring homomorphisms and ring homomorphisms of finite flat dimension. Both of these types of
ring homomorphisms contain the class of faithfully flat ring homomorphisms. In fact, by \cite[Theorem 7.5]{M}, an
$R$-algebra $T$ is faithfully flat if and only if it is both pure and flat. We proceed with the paper as follows.

In Section 2, we summarize some needed definitions and results from \cite{DGTZ2} and \cite{PDR}.

Section 3 deals with pure ring homomorphisms. Let $f:R\lo T$ be a pure ring homomorphism. We show that if the $T$-module
$T\otimes_RM$ is $\fa T$-relative Cohen-Macaulay (resp. $\fa T$-relative maximal Cohen-Macaulay), then $M$ is
$\fa$-relative Cohen-Macaulay (resp. $\fa$-relative maximal Cohen-Macaulay). Under some extra assumptions, we also prove
that if $T\otimes_RM$ is $\fa T$-relative regular, then $M$ is $\fa$-relative regular; see Theorem \ref{3.3}.  Example
\ref{3.14} indicates that Theorem \ref{3.3} fails for relative Gorensteinness. Also, we illustrate that Theorem \ref{3.3}
fails if $f$ has finite flat dimension instead of being pure; see Example \ref{3.4}.

Section 4 discusses ring homomorphisms of finite flat dimension. Let $f:R\lo T$ be a flat ring homomorphism with
$\fa T\neq T$. We prove that if $M$ is $\fa$-relative Cohen-Macaulay (resp. $\fa$-relative Gorenstein,  $\fa$-relative
regular), then $T\otimes_RM$ is $\fa T$-relative Cohen-Macaulay (resp. $\fa T$-relative Gorenstein, $\fa T$-relative
regular). Also, we show that if $M$ is $\fa$-relative maximal Cohen-Macaulay and $T\otimes_RM\neq (\fa T)(T\otimes_RM)$,
then $T\otimes_RM$ is $\fa T$-relative maximal Cohen-Macaulay, see Theorem \ref{4.2}. Theorem \ref{4.2} fails if we
replace the assumption that $f$ is flat with the weaker assumption that $f$ has finite flat dimension; see Example
\ref{4.4}. Theorems \ref{3.3} and \ref{4.2} allow us to derive Proposition \ref{4.6} concerning faithfully flat ring
homomorphisms:  Let $f:R\lo T$ be a faithfully flat ring homomorphism. Then, $M$ is $\fa$-relative Cohen-Macaulay
(resp. $\fa$-relative maximal Cohen-Macaulay, $\fa$-relative Gorenstein) if and only if $T\otimes_RM$ is $\fa T$-relative
Cohen-Macaulay (resp. $\fa T$-relative maximal Cohen-Macaulay, $\fa T$-relative Gorenstein). Section 4 continues with
considering the ascent and descent  of relative homological modules along the natural ring epimorphism
$\pi:R\lo R/\langle x_1,x_2, ..., x_n\rangle$, where $x_1,x_2, ..., x_n$ is an $R$-regular sequence.

\section{Prerequisites}

The purpose of this section is to provide some background material on relative systems of parameters and $\fa$-relative
Cohen-Macaulay, $\fa$-relative maximal Cohen-Macaulay, $\fa$-relative Gorenstein and $\fa$-relative regular modules, which
will be used in the rest of the paper.

Below, we recall some definitions and results from \cite{DGTZ2}. Let $N$ be an $R$-module. This paper addresses, among other
things, local cohomology modules $$\RH_{\fa}^{i}\left(N\right)=\varinjlim \limits_{n\in \mathbb{N}}
\text{Ext}_R^i\left(R/\fa^n,N\right); \  i\in \mathbb{N}_0.$$ We denote by $\cd(\fa,N)$ the {\it cohomological dimension of
$N$ with respect to $\fa$} which is the largest integer $i$ such that $\RH_{\fa}^{i}\left(N\right)\neq 0$. One may easily check
that $\cd\left(\fa,N\right)\leq \ara\left(\fa(R/\Ann_RN)\right)$ and $\cd\left(\fa,N\right)\leq \cd\left(\fa,R\right)$.
(The arithmetic rank of an ideal $\fb$ of a ring $T$ is denoted by $\ara(\fb)$.)

\begin{definition}\label{2.1} Let $M$ be a finitely generated $R$-module and $\fa$ an ideal of $R$ with $M\neq \fa M$.
Let $c=\cd\left(\fa,M\right)$. A sequence $x_{1}, x_2, \ldots, x_{c}\in \fa$ is called $\fa$-{\it relative system of
parameters}, $\fa$-s.o.p, of $M$ if $$\Rad\left(\langle x_{1}, x_2, \ldots, x_{c}\rangle+\Ann_{R}M\right)=
\Rad\left(\fa+\Ann_{R}M\right).$$
\end{definition}

Clearly, if $x_{1},x_2,\ldots, x_{c}\in R$ is an $\fa$-s.o.p of $M$, then for any $t_{1},\ldots,t_{c}\in \mathbb{N}$, every
permutation of $x_{1}^{t_{1}},\ldots, x_{c}^{t_{c}}$ is also an $\fa$-s.o.p of $M$.

Over a local ring, every finitely generated module possesses a system of parameters, but this is not true for $\fa$-relative
systems of parameters. In fact, by \cite[Lemma 2.2]{DGTZ2}, for a finitely generated $R$-module $M$ and an ideal $\fa$ of $R$
with $M\neq \fa M$, $\fa$ contains an $\fa$-s.o.p of $M$ if and only if $\cd\left(\fa,M\right)=\ara\left(\fa(R/\Ann_RM)\right)$.

\begin{theorem}\label{2.2} (See \cite[Lemma 2.4 and Theorem 2.7]{DGTZ2}.) Let $\fa$ be an ideal of $R$, $M$ a finitely generated
$R$-module with $M\neq \fa M$ and $c=\cd\left(\fa,M\right)$. Assume that $\cd\left(\fa,M\right)=\ara\left(\fa(R/\Ann_RM)\right)$
and $x_1,\ldots, x_{c}\in \fa$. Consider the following conditions:
\begin{enumerate}
\item[(i)] $x_1,\ldots, x_{c}$ is an $\fa$-s.o.p of $M$.
\item[(ii)] $\cd\left(\fa,M/\langle x_{1},x_2, \ldots, x_{i}\rangle M\right)=c-i$ for every $i=1, 2,\ldots, c$.
\end{enumerate}
Then (i) implies (ii). Furthermore, if $\fa$ is contained in the Jacobson radical of $R$, then (i) and (ii) are equivalent.
\end{theorem}

\begin{definition}\label{2.3} Let $\fa$ be a proper ideal of $R$.
\begin{enumerate}
\item[(i)] A finitely generated $R$-module $M$ is called $\fa$-{\it relative Cohen-Macaulay} if either $M=\fa M$ or $M\neq \fa M$
and $\grade(\fa,M)=\cd(\fa,M)$.
\item[(ii)] A finitely generated $R$-module $M$ is called $\fa$-{\it relative maximal Cohen-Macaulay} if $\grade(\fa,M)=\cd(\fa,R)$.
\item[(iii)] A finitely generated $R$-module $M$ is called $\fa$-{\it relative Gorenstein} if $\Ext_R^i(R/\fa,M)=0$ for all
$i\neq \cd(\fa,R)$.
\item[(iv)] A finitely generated $R$-module $M$ is called $\fa$-{\it relative regular} if either $M=\fa M$ or $M\neq \fa M$ and
$\fa$ can be generated by an $M$-regular sequence of length $\cd(\fa,R)$ which is also an $R$-regular sequence.
\end{enumerate}
\end{definition}

Let $(R,\fm)$ be a local ring and $M$ a finitely generated $R$-module. Then $R$ is regular if and only if $R$ is $\fm$-relative
regular. Also, $M$ is Gorenstein (resp. maximal Cohen-Macaulay, Cohen-Macaulay) if and only if $M$ is $\fm$-relative Gorenstein
(resp. $\fm$-relative maximal Cohen-Macaulay, $\fm$-relative Cohen-Macaulay).

\begin{theorem}\label{2.4} (See \cite[Theorem 2.19 and Proposition 4.6]{PDR}.) Let $\fa$ be a proper ideal of $R$.
\begin{enumerate}
\item[(i)]  Every $\fa$-relative regular module is $\fa$-relative Gorenstein.
\item[(ii)]  Every $\fa$-relative Gorenstein module $M$ with $M\neq \fa M$ is $\fa$-relative maximal Cohen-Macaulay.
\item[(iii)] Every $\fa$-relative maximal Cohen-Macaulay module is $\fa$-relative Cohen-Macaulay.
\end{enumerate}
\end{theorem}

\begin{theorem}\label{2.5} (See \cite[Theorem 3.3]{DGTZ2}.) Let $M$ be a finitely generated $R$-module and $\fa$ a proper ideal
of $R$ with $\cd\left(\fa,M\right)=\ara\left(\fa(R/\Ann_RM)\right)$. Consider the following conditions:
\begin{enumerate}
\item[(i)] $M$ is $\fa$-relative Cohen-Macaulay.
\item[(ii)] Every $\fa$-s.o.p of $M$ is an $M$-regular sequence.
\item[(iii)] There exists an $\fa$-s.o.p of $M$ which is an $M$-regular sequence.
\end{enumerate}
Then (i) and (iii) are equivalent. Furthermore, if $\fa$ is contained in the Jacobson radical of $R$, then all three conditions
are equivalent.
\end{theorem}

\section{Pure ring homomorphisms}

In this section, we will investigate the behaviour of relative homological modules under pure ring homomorphisms.

\begin{definition}\label{3.1} A ring homomorphism $f:R\lo T$ is called {\it pure} if the $R$-homomorphism $$f\otimes
\id_M:R\otimes_RM\lo T\otimes_RM$$ is injective for every $R$-module $M$.
\end{definition}

Clearly, every pure ring homomorphism is injective.

\begin{lemma}\label{3.2}  Let $f:R\lo T$ be a pure ring homomorphism, $\fa$ a proper ideal of $R$ and $M$ a finitely generated
$R$-module. Then
\begin{enumerate}
\item[(i)] $\grade(\fa,M)\geq \grade(\fa T,T\otimes_RM)$.
\item[(ii)] $\cd\left(\fa,M\right)=\cd\left(\fa T,T\otimes_RM\right)$.
\end{enumerate}
\end{lemma}

\begin{prf} For any non-negative integer $i$, \cite[Corollary 6.8]{HR} asserts that the induced $R$-homomorphism
$$\text{H}_{\fa}^i(f\otimes \id_M): \text{H}_{\fa}^i(R\otimes_RM)\lo \text{H}_{\fa T}^i(T\otimes_RM)$$ is injective.
Thus, $$\{i\in \mathbb{N}_0 \mid \text{H}_{\fa}^i(M)\neq 0\}\subseteq \{i\in \mathbb{N}_0 \mid
\text{H}_{\fa T}^i(T\otimes_RM)\neq 0\}.$$ Hence, $\grade(\fa T,T\otimes_RM)\leq \grade(\fa,M)$ and $\cd(\fa,M)\leq
\cd(\fa T,T\otimes_RM).$

Recall that local cohomology functors commute with direct limits, and every $R$-module is the direct limit of its finitely
generated submodules. Thus, by \cite[Theorem 2.2]{DNT}, we conclude that if $X$ and $Y$ are two $R$-modules such that $Y$
is finitely generated and $\Supp_RX\subseteq \Supp_RY$, then $\cd(\fa,X)\leq \cd(\fa,Y)$. The Independent theorem
\cite[Theorem 4.2.1]{BS} yields an $R$-isomorphism $\text{H}_{\fa T}^i(N)\cong \text{H}_{\fa}^i(N)$ for every $T$-module
$N$ and all $i\geq 0$. Thus, $$\cd(\fa T,T\otimes_RM)=\cd(\fa,T\otimes_RM)\leq \cd(\fa,M),$$ and so $\cd(\fa T,T\otimes_RM)
=\cd(\fa,M)$.
\end{prf}

The next result might be considered as a relative analogue of the Cohen-Macaulay Direct Summand theorem; see e.g.
\cite[Theorem 10.4.1]{BH}.

In what follows, $\mu(\fa)$ stands for the minimum number of generators of an ideal $\fa$.

\begin{theorem}\label{3.3}  Let $f:R\lo T$ be a pure ring homomorphism, $\fa$ a proper ideal of $R$ and $M$ a finitely
generated $R$-module.
\begin{enumerate}
\item[(i)]  If $T\otimes_RM$ is $\fa T$-relative Cohen-Macaulay, then $M$ is $\fa$-relative Cohen-Macaulay.
\item[(ii)] If $T\otimes_RM$ is $\fa T$-relative maximal Cohen-Macaulay, then $M$ is $\fa$-relative maximal Cohen-Macaulay.
\item[(iii)] Assume that $\fa$ is contained in the Jacobson radical of $R$ and $\mu(\fa)=\cd(\fa,R)$. If $T\otimes_RM$ is
$\fa T$-relative regular, then $M$ is $\fa$-relative regular.
\end{enumerate}
\end{theorem}

\begin{prf} As $f$ is pure, an $R$-module $N$ is zero if and only if the $T$-module $T\otimes_RN$ is zero. In particular,
$\fa T$ is a proper ideal of $T$ and $M=\fa M$ if and only if $T\otimes_RM=(\fa T)(T\otimes_RM)$.

(i) Suppose that $T\otimes_RM$ is $\fa T$-relative Cohen-Macaulay. If $M=\fa M$, then $M$ is $\fa$-relative Cohen-Macaulay
by the definition. So, we may assume that $M\neq \fa M$.  Then, $T\otimes_RM\neq (\fa T)(T\otimes_RM)$. Next, Lemma \ref{3.2}
yields the following display:
$$\begin{array}{ll}
\cd(\fa T,T\otimes_RM)&=\grade(\fa T,T\otimes_RM)\\
&\leq \grade(\fa,M)\\
&\leq \cd(\fa,M)\\
&=\cd(\fa T,T\otimes_RM).
\end{array}$$
Thus $\grade(\fa,M)=\cd(\fa,M)$, and so $M$ is $\fa$-relative Cohen-Macaulay.

(ii)  Suppose that $T\otimes_RM$ is $\fa T$-relative maximal Cohen-Macaulay. In particular, this yields that $T\otimes_RM\neq
(\fa T)(T\otimes_RM)$, and so $M\neq \fa M$. Now, Lemma \ref{3.2} implies the following display:
$$\begin{array}{ll}
\cd(\fa T,T)&=\grade(\fa T,T\otimes_RM)\\
&\leq \grade(\fa,M)\\
&\leq \cd(\fa,M)\\
&\leq \cd(\fa,R)\\
&=\cd(\fa T,T).
\end{array}$$
Hence $\grade(\fa,M)=\cd(\fa,R)$, and so $M$ is $\fa$-relative maximal Cohen-Macaulay.

(iii) Suppose that $T\otimes_RM$ is $\fa T$-relative regular, and set $c=\cd(\fa,R)$. If $T\otimes_RM=(\fa T)(T\otimes_RM)$,
then $M=\fa M$, and so $M$ is $\fa$-relative regular. Hence, we may and do assume that $T\otimes_RM\neq (\fa T)(T\otimes_RM)$.
Now, by the definition, it follows that $T$ is $\fa T$-relative regular. Theorem \ref{2.4} implies that $T\otimes_RM$ is
$\fa T$-relative maximal Cohen-Macaulay. Thus, by (ii), we conclude that $M$ is $\fa$-relative maximal Cohen-Macaulay. As a
result, $\cd(\fa,M)=c$. Similarly, by Theorem \ref{2.4} and (i), we deduce that $R$ is $\fa$-relative Cohen-Macaulay.

Since $\mu(\fa)=c$, there are elements $x_1,..., x_c\in R$ such that $\fa=\langle x_1,..., x_c \rangle$.  It is evident that
$x_1,...,x_c$ is an $\fa$-s.o.p of both $M$ and $R$. Now, as $\fa$ is contained in the Jacobson radical of $R$, Theorem
\ref{2.5} implies that $x_1,..., x_c$ is both an $M$-regular sequence and an $R$-regular sequence. Therefore, $M$ is
$\fa$-relative regular.
\end{prf}

Example \ref{3.14} illustrates Theorem \ref{3.3}'s failure for relative Gorensteinness.

Let $\fa$ be an ideal of $R$ contained in its Jacobson radical. If $R$ admits a nonzero $\fa$-relative regular module, then
\cite[Theorem 4.4]{PDR} yields that $\mu(\fa)=\cd(\fa,R)$. So, the assumption $\mu(\fa)=\cd(\fa,R)$ in Theorem \ref{3.3}(iii)
is justified.

In Theorem \ref{3.3}, one may ask whether it is possible to replace the assumption that $f$ is pure with the that $f$ has
finite flat dimension. This is not the case, as shown in the following example.

\begin{example}\label{3.4}  Let $R=\Bbbk[x,y]$ be a polynomial ring over a field $\Bbbk$ and set $T=R/\langle y\rangle$.
Let $f:R\to T$ be the natural ring epimorphism. As $\pd_RT=1$, we have $\fd_RT<\infty$. We set $\fa=\langle x \rangle$ and
$M=R/\langle xy\rangle$. Observes $T\otimes_RM\cong T.$ Since $\fa T=\langle x+ \langle y\rangle \rangle$ and $x+ \langle
y\rangle$ is $T$-regular, it follows that $\cd(\fa T, T)=1$, and so $T\otimes_RM$ is $\fa T$-relative regular. Now, of course,
$T\otimes_RM$ is $\fa T$-relative Gorenstein, $\fa T$-relative maximal Cohen-Macaulay and $\fa T$-relative Cohen-Macaulay.

On the other hand,  we claim that $M$ is not $\fa$-relative Cohen-Macaulay. As $\Ass_RM=\{\langle x\rangle,\langle y\rangle\}$,
one has $$\cd(\fa,M)=\max\{\cd(\fa,R/\langle x\rangle),\cd(\fa,R/\langle y\rangle)\}=\max\{0,1\}=1.$$ Since $\fa=\langle x\rangle
\in \Ass_RM$, it follows that $\grade(\fa,M)=0$. Hence, $M$ is not $\fa$-relative Cohen-Macaulay. Thus, $M$ is not $\fa$-relative
maximal Cohen-Macaulay, not $\fa$-relative Gorenstein and not $\fa$-relative regular.
\end{example}

\begin{lemma}\label{3.5}  Let $\fa$ and $\fb$ be two proper ideals of $R$ with the same radical and $M$ a finitely generated
$R$-module. Then
\begin{enumerate}
\item[(i)] $M$ is $\fa$-relative Cohen-Macaulay if and only if $M$ is $\fb$-relative Cohen-Macaulay.
\item[(ii)] $M$ is $\fa$-relative maximal Cohen-Macaulay if and only if $M$ is $\fb$-relative maximal Cohen-Macaulay.
\end{enumerate}
\end{lemma}

\begin{prf} As $\Rad(\fa)=\Rad(\fb)$, it follows that $\grade(\fa,M)=\grade(\fb,M)$, $\cd(\fa,M)=\cd(\fb,M)$ and $\cd(\fa,R)=
\cd(\fb,R)$. These observations conclude the assertions.
\end{prf}

The following immediate corollary, concerning homological modules, is interesting in its own right.

\begin{corollary}\label{3.6}  Let $f:(R,\fm)\lo (T,\fn)$ be a pure ring homomorphism of local rings and $M$ a finitely generated
$R$-module. Assume that $\Rad(\fm T)=\fn$. Then $\dim_T(T\otimes_RM)=\dim_RM$. Furthermore:
\begin{enumerate}
\item[(i)]  If $T\otimes_RM$ is a Cohen-Macaulay $T$-module, then $M$ is a Cohen-Macaulay $R$-module.
\item[(ii)] If $T\otimes_RM$ is a maximal Cohen-Macaulay $T$-module, then $M$ is a maximal Cohen-Macaulay $R$-module.
\end{enumerate}
\end{corollary}

\begin{prf} By Grothendieck's Non-vanishing theorem \cite[Theorem 6.1.4]{BS}, if $N$ is a finitely generated module over a
local ring $(S,\fl)$, then $\cd(\fl,N)=\dim_SN$. Hence, $\dim_T(T\otimes_RM)=\dim_RM$ by Lemma \ref{3.2}(ii). The rest of
the claim follows from Theorem \ref{3.3} and Lemma \ref{3.5}.
\end{prf}

Let $\fa$ and $\fb$ be two proper ideals of $R$ with the same radical and $M$ a finitely generated $R$-module. Keeping Lemma
\ref{3.5} in mind, one may guess that $M$ is $\fa$-relative Gorenstein (resp. $\fa$-relative regular) if and only if $M$ is
$\fb$-relative Gorenstein (resp. $\fb$-relative regular). As the following example illustrates, this is not the case.

\begin{example}\label{3.7} Let $(R,\fm)$ be a non-Gorenstein Cohen-Macaulay local ring of dimension $d$. Let $x_1, x_2, \dots,
x_d\in \fm$ be a system of parameters of $R$, and set $\fa=\langle x_1, x_2, \dots, x_d\rangle$. Then $\Rad(\fa)=\fm$, and so
$$\cd(\fa,R)=\cd(\fm,R)=\dim R=d.$$ As $R$ is Cohen-Macaulay, it follows that $x_1,x_2, \dots, x_d$ is an $R$-regular sequence.
Hence, $R$ is $\fa$-relative regular, and so it is $\fa$-relative Gorenstein as well. However, $R$ is not Gorenstein, and therefore
it is not regular.
\end{example}

Below, we recall the definitions of ring retractions and trivial extensions.

\begin{definition}\label{3.8} A ring monomorphism $f:R\lo T$ is said to be a {\it ring retraction} if there is a ring homomorphism
$g:T\lo R$ such that $gf=\id_R$.
\end{definition}

\begin{definition}\label{3.9} Let $M$ be an $R$-module. Then $R\oplus M$ with coordinate-wise addition and multiplication $(r_1,m_1)(r_2,m_2)=(r_1r_2,r_1m_2+r_2m_1)$ is a commutative ring with identity. It is called {\it trivial extension} of
$R$ by $M$ and denoted by $R\ltimes M$.
\end{definition}

If $M$ is a finitely generated $R$-module, then the ring $R\ltimes M$ is Noetherian by \cite[Proposition 2.2]{AW}.

\begin{remark}\label{3.10}
\begin{enumerate}
\item[(i)] Let $f:R\lo T$ be a ring homomorphism. If either $f$ is faithfully flat, $f$ is a ring retraction, or $R$ is a direct
summand of $T$ as an $R$-module, then $f$ is pure. Thus, the assertions of Lemma \ref{3.2} and Theorem \ref{3.3} hold for $f$,
in each of these three cases.
\item[(ii)] Let $M$ be a finitely generated $R$-module, and set $T=R\ltimes M$. It is easy to verify that the canonical
injection  $\lambda:R\lo T$ and the canonical projection $\rho: T\to R$ are ring homomorphism, and $\rho\lambda=\id_R$.
Hence, $\lambda: R\lo T$ is a ring retraction, and so it is pure by (i). Thus, the assertions of Lemma \ref{3.2} and Theorem
\ref{3.3} hold for $\lambda$.
\item[(iii)]  Let $R=\Bbbk[x_1,x_2,...,x_n]$ be a polynomial ring over a field $\Bbbk$ and $G$ a finite group acting on $R$.
Let $R^G$ be the invariant subring of $R$ under this action. Assume that the characteristic of $\Bbbk$ is either zero or a
prime number not dividing $|G|$. Consider Reynold's operation $\theta: R\lo R^G$ given by $\theta(x)=\frac{1}{|G|}\underset{g\in
G}\sum g.x$ and let $\iota: R^G\lo R$ denote the natural inclusion homomorphism. Then $\theta\iota=\id_{R^G}$, and so
$\iota: R^G\lo R$ is a ring retraction.
\end{enumerate}
\end{remark}

Next, we improve Lemma \ref{3.2} and Theorem \ref{3.3} for trivial extensions of $R$.

\begin{lemma}\label{3.11}  Let $\fa$ be an ideal of $R$, $M$ a finitely generated $R$-module and let $T=R\ltimes M$ denote the
trivial extension of $R$ by $M$. Let $\fa T$ be the extension of the ideal $\fa$ under the canonical injection $\lambda: R\to T$.
Then
\begin{enumerate}
\item[(i)] $\grade(\fa T,T)=\min \{\grade(\fa, R),\grade(\fa,M) \}$.
\item[(ii)]  $\cd(\fa T,T)=\cd(\fa,R)$.
\end{enumerate}
\end{lemma}

\begin{prf} Let $\lambda: R\lo T$ and $\rho: T\to R$ be, respectively,  the canonical injection and canonical projection.
Then, every $T$-module admits the structure of an $R$-module via $\lambda$ and every $R$-module admits the structure of a
$T$-module via $\rho$. For an $R$-module $L$, if we equip it with the structure of a $T$-module via $\rho$, and then equip the
resulting $T$-module with the structure of an $R$-module via $\lambda$, we regain the original $R$-module $L$
that we started with.

(i) Since as an $R$-module, we have $T=R\oplus M$, it follows that $$\grade(\fa,T)=\min \{\grade(\fa,R),\grade(\fa,M) \}.$$ For
every finitely generated $T$-module $N$, by the Independent theorem, we deduce that $\grade(\fa T,N)=\grade(\fa,N)$, and so (i)
follows. Note that in light of $T$ being finitely generated as an $R$-module, every finitely generated $T$-module is also finitely
generated when it is viewed as an $R$-module via $\lambda$.

(ii) is immediate by Remark \ref{3.10}(ii) and Lemma \ref{3.2}(ii).
\end{prf}

\begin{proposition}\label{3.12} Let $\fa$ be an ideal of $R$, $M$ a finitely generated $R$-module with $M\neq \fa M$ and let
$T=R\ltimes M$ denote the trivial extension of $R$ by $M$.
\begin{enumerate}
\item[(i)] $T$ is $\fa T$-relative Cohen-Macaulay if and only if $R$ is $\fa$-relative Cohen-Macaulay and $M$ is $\fa$-relative
maximal Cohen-Macaulay.
\item[(ii)] Assume that $\fa$ is contained in the Jacobson radical of $R$ and $\mu(\fa)=\cd(\fa,R)$. Then $T$ is $\fa T$-relative
regular if and only if $R$ is $\fa$-relative regular and $M$ is $\fa$-relative maximal Cohen-Macaulay if and only if $M$ is
$\fa$-relative regular.
\end{enumerate}
\end{proposition}

\begin{prf} Note that, by Remark \ref{3.10}(ii), the canonical injection $\lambda: R\lo T$ is pure. Also, note that $\cd(\fa T,T)=
\cd(\fa,R)$ by Lemma \ref{3.11}(ii).

(i) Assume that $T$ is $\fa T$-relative Cohen-Macaulay. As $\lambda$ is pure, Theorem \ref{3.3}(i) asserts that $R$ is
$\fa$-relative Cohen-Macaulay. Next, Lemma \ref{3.11}(i) yields
$$\begin{array}{ll}
\cd(\fa T,T)&=\grade(\fa T,T)\\
&\leq \grade(\fa,M)\\
&\leq \cd(\fa,M)\\
&\leq \cd(\fa,R)\\
&=\cd(\fa T,T).
\end{array}$$
This implies that $\grade(\fa,M)=\cd(\fa,R)$, and so $M$ is $\fa$-relative maximal Cohen-Macaulay.

Conversely, assume that $R$ is $\fa$-relative Cohen-Macaulay and $M$ is $\fa$-relative maximal Cohen-Macaulay. Then $$\grade(\fa,M)
=\cd(\fa,R)=\grade(\fa,R).$$ Now, Lemma \ref{3.11}(i) implies that $$\grade(\fa T,T)=\grade(\fa,R)=\cd(\fa,R)=\cd(\fa T,T),$$
and so $T$ is $\fa T$-relative Cohen-Macaulay.

(ii) Set $c=\cd(\fa,R)$.  Assume that $R$ is $\fa$-relative regular and $M$ is $\fa$-relative maximal Cohen-Macaulay. Then
\cite[Lemma 4.2]{PDR} implies that $M$ is $\fa$-relative regular. Hence, $\fa$ can be generated by some elements $x_1, \dots,
x_c$ which form both an $M$-regular sequence and an $R$-regular sequence. Now, it is straightforward to check that $(x_1,0),
\dots, (x_c,0)$ is a $T$-regular sequence and $\fa T=\langle (x_1,0), \dots, (x_c,0) \rangle$. Thus, $T$ is $\fa T$-relative
regular.

Now, assume that $T$ is $\fa T$-relative regular. Then, $R$ is $\fa$-relative regular by Theorem \ref{3.3}(iii). As $T$
is $\fa T$-relative regular, Theorem \ref{2.4} implies that $T$ is $\fa T$-relative Cohen-Macaulay. Now, (i) yields that $M$ is
$\fa$-relative maximal Cohen-Macaulay.

As $M\neq \fa M$, the $\fa$-relative regularness of $M$ implies that $R$ is $\fa$-relative regular. So, by \cite[Lemma 4.2]{PDR}
and Theorem \ref{2.4}, it turns out that $R$ is $\fa$-relative regular and $M$ is $\fa$-relative maximal Cohen-Macaulay if and only
if $M$ is $\fa$-relative regular.
\end{prf}

Let $f: R\lo T$ be a pure ring homomorphism, $\fa$ a proper ideal of $R$ and $M$ a finitely generated $R$-module. In light of
Theorem \ref{3.3}, one may conjecture that if $T\otimes_RM$ is $\fa T$-relative Gorenstein, then $M$ is $\fa$-relative Gorenstein.
This section concludes with an example that shows this is not the case. It requires the following lemma.

Recall that a finitely generated $R$-module $C$ is called {\it semidualizing} if the homothety map $\chi_C^R:R \lo
\Hom_R\left(C,C\right)$ is an isomorphism, and $\Ext^i_R\left(C,C\right)=0$ for all $i>0$.

\begin{lemma}\label{3.13}  Let $\fa$ be a proper ideal of $R$, $C$ a semidualizing module of $R$ and let $T=R\ltimes C$ denote the
trivial extension of $R$ by $C$. Let $\fa T$ be the extension of the ideal $\fa$ under the canonical injection $\lambda: R\to T$.
Then $T$ is $\fa T$-relative Gorenstein if and only if $C$ is $\fa$-relative Gorenstein and $\Ext_R^i(C/\fa C,C)=0$ for all $i\neq
\cd(\fa,R)$.
\end{lemma}

\begin{prf} By \cite[Lemma 3.2(ii)]{HJ}, we get an $R$-isomorphism $\Ext_T^i(L,T)\cong \Ext_R^i(L,C)$ for every $R$-module $L$ and all
$i\geq 0$. Because of the $R$-isomorphism $T/\fa T\cong R/\fa\oplus C/\fa C$, we deduce an $R$-isomorphism $$\Ext_T^i(T/\fa T,T)
\cong \Ext_R^i(R/\fa,C)\oplus \Ext_R^i(C/\fa C,C)$$ for all $i\geq 0$. In particular, for each integer $i\geq 0$, $\Ext_T^i(T/\fa T,T)$
vanishes if and only if both $\Ext_R^i(R/\fa,C)$ and $\Ext_R^i(C/\fa C,C)$ vanish. This completes the proof, because $\cd(\fa T,T)=
\cd(\fa, R)$ by Lemma \ref{3.11}(ii).
\end{prf}

\begin{example}\label{3.14} Let $(R,\fm,\Bbbk)$ be a non-Gorenstein Cohen-Macaulay local ring with a dualizing module $\omega_R$. Then,
$\omega_R$ is a semidualizing module of $R$ and it is $\fm$-relative Gorenstein. Set $d=\dim R=\cd(\fm,R)$. Hence, $\Ext_R^i(\Bbbk,
\omega_R)=0$ for all $i\neq d$. As $\omega_R/\fm \omega_R$ is isomorphic, as an $R$-module, to a finite direct sum of copies of $\Bbbk$,
it turns out that $\Ext_R^i(\omega_R/\fm \omega_R,\omega_R)=0$ for all $i\neq d$. Now, by Lemma \ref{3.13}, $T$ is $\fm T$-relative
Gorenstein, while $R$ is not $\fm$-relative Gorenstein.
\end{example}

\section{Ring homomorphisms of finite flat dimension}

In this section, we will examine the behaviour of relative homological modules under ring homomorphisms of finite flat dimension.

\begin{lemma}\label{4.1} Let $f:R\lo T$ be a flat ring homomorphism and $\fa$ a proper ideal of $R$. Let $M$ be a finitely
generated $R$-module with $T\otimes_RM\neq  (\fa T)(T\otimes_RM)$.
\begin{enumerate}
\item[(i)] If $M$ is $\fa$-relative Cohen-Macaulay, then $\cd\left(\fa T,T\otimes_RM\right)=\cd\left(\fa,M \right)$.
\item[(ii)] If $M$ is $\fa$-relative maximal Cohen-Macaulay, then $\cd\left(\fa T,T\right)=\cd\left(\fa,R \right)$.
\end{enumerate}
\end{lemma}

\begin{prf} The Flat Base Change theorem \cite[Theorem 4.3.2]{BS} yields a $T$-isomorphism
$${\RH}_{\fa T}^{i}\left(T\otimes_RM\right)\cong {\RH}_{\fa}^{i}\left(M\right)\otimes_RT \   \  (*)$$ for all $i\geq 0$.
As $T\otimes_RM\neq  (\fa T)(T\otimes_RM)$, obviously, we have $M\neq \fa M$. Hence, $\cd\left(\fa T,T\otimes_RM\right)$ and
$\cd\left(\fa,M \right)$ are non-negative integers.

(i) Suppose that $M$ is $\fa$-relative Cohen-Macaulay. Then, by $(*)$, we conclude that $\cd\left(\fa T,T\otimes_RM\right)=
\cd\left(\fa,M \right)$.

(ii) Suppose that $M$ is $\fa$-relative maximal Cohen-Macaulay. Then, $M$ is $\fa$-relative Cohen-Macaulay by Theorem \ref{2.4}(iii).
Now, by (i), one has the following display:
$$\begin{array}{ll}
\cd(\fa,R)&=\grade(\fa,M)\\
&=\cd(\fa,M)\\
&=\cd(\fa T,T\otimes_RM)\\
&\leq \cd(\fa T,T)\\
&=\cd(\fa,T)\\
&\leq \cd(\fa,R).
\end{array}$$
Hence, $\cd\left(\fa T,T\right)=\cd\left(\fa,R \right)$.
\end{prf}

\begin{theorem}\label{4.2} Let $f:R\lo T$ be a flat ring homomorphism, $M$ a finitely generated $R$-module and $\fa$ an
ideal of $R$ with $\fa T\neq T$.
\begin{enumerate}
\item[(i)] If $M$ is $\fa$-relative Cohen-Macaulay, then $T\otimes_RM$ is $\fa T$-relative Cohen-Macaulay.
\item[(ii)] If $M$ is $\fa$-relative maximal Cohen-Macaulay and $T\otimes_RM\neq (\fa T)(T\otimes_RM)$, then $T\otimes_RM$ is
$\fa T$-relative maximal Cohen-Macaulay.
\item[(iii)] If $M$ is $\fa$-relative Gorenstein, then $T\otimes_RM$ is $\fa T$-relative Gorenstein.
\item[(iv)] If $M$ is $\fa$-relative regular, then $T\otimes_RM$ is $\fa T$-relative regular.
\end{enumerate}
\end{theorem}

\begin{prf} (i) Assume that $M$ is $\fa$-relative Cohen-Macaulay. If $T\otimes_RM=(\fa T)(T\otimes_RM)$, then $T\otimes_RM$
is $\fa T$-relative Cohen-Macaulay by the definition. So, we may assume that $T\otimes_RM\neq (\fa T)(T\otimes_RM)$. Then the
isomorphisms $(*)$, in the proof of Lemma \ref{4.1}, yields that $$\grade(\fa,M)\leq \grade(\fa T,T\otimes_RM).$$ Now, by
Lemma \ref{4.1}(i), we have the following display:
$$\begin{array}{ll}
\cd(\fa T,T\otimes_RM)&=\cd(\fa,M)\\
&=\grade(\fa,M)\\
&\leq \grade(\fa T,T\otimes_RM)\\
&\leq \cd(\fa T,T\otimes_RM).
\end{array}$$
Hence, $$\grade(\fa T,T\otimes_RM)=\cd(\fa T,T\otimes_RM),$$ and so $T\otimes_RM$ is $\fa T$-relative Cohen-Macaulay.

(ii) Assume that $M$ is $\fa$-relative maximal Cohen-Macaulay and $T\otimes_RM\neq  (\fa T)(T\otimes_RM)$. By Lemma \ref{4.1}(ii),
we deduce the following display:
$$\begin{array}{ll}
\cd(\fa T,T)&= \cd(\fa,R)\\
&=\grade(\fa,M)\\
&\leq \grade(\fa T,T\otimes_RM)\\
&\leq \cd(\fa T,T\otimes_RM)\\
&\leq \cd(\fa T,T).
\end{array}$$
Thus, $$\grade(\fa T,T\otimes_RM)=\cd(\fa T,T),$$ and so $T\otimes_RM$ is $\fa T$-relative maximal Cohen-Macaulay.

(iii) Assume that $M$ is $\fa$-relative Gorenstein. If $T\otimes_RM=(\fa T)(T\otimes_RM)$, then $$\Ext_T^i(T/\fa T,T\otimes_RM)=0$$
for all $i\geq 0$, and so $T\otimes_RM$ is $\fa T$-relative Gorenstein. So, we may assume that $T\otimes_RM\neq  (\fa T)(T\otimes_RM)$.
This, in particular, yields that $M\neq \fa M$, and so $M$ is $\fa$-relative maximal Cohen-Macaulay by Theorem \ref{2.4}(ii). Now,
Lemma \ref{4.1}(ii) implies that $\cd(\fa T,T)=\cd(\fa,R)$. As $f$ is flat, there is a natural $T$-isomorphism $$\Ext_T^i(T/\fa T,T\otimes_RM)\cong \Ext_R^i(R/\fa,M)\otimes_RT$$ for all $i\geq 0$. Since, $M$ is $\fa$-relative Gorenstein, it turns out that
$\Ext_T^i(T/\fa T,T\otimes_RM)=0$ for all $i\neq \cd\left(\fa T,T\right)$, and so $T\otimes_RM$ is $\fa T$-relative Gorenstein.

(iv) Assume that $M$ is $\fa$-relative regular. If $T\otimes_RM=(\fa T)(T\otimes_RM)$, then $T\otimes_RM$ is $\fa T$-relative
regular by the definition. Hence, we may and do assume that $T\otimes_RM\neq (\fa T)(T\otimes_RM)$, and so $M\neq \fa M$. Let
$c=\cd\left(\fa,R\right)$. By our assumption, $\fa$ has generators $x_{1}, x_2, \dots, x_{c}$ which form both an $M$-regular 
sequence and an $R$-regular sequence. As $\fa=\langle x_{1}, x_2, \ldots, x_{c}\rangle$, one has $\fa T=\langle f(x_{1}), f(x_2), \ldots,
f(x_{c})\rangle$. Since $f$ is flat, it is routine to verify that $f(x_{1}), f(x_2), \ldots, f(x_{c})$ is both a
$T\otimes_RM$-regular sequence and a $T$-regular sequence. Thus, $T\otimes_RM$ is $\fa T$-relative regular. Note that, by Theorem
\ref{2.4} and Lemma \ref{4.1}(ii), $\cd(\fa T,T)=c$.
\end{prf}

Next, we record the following immediate corollary.

\begin{corollary}\label{4.3} Let $\fa$ be an ideal of $R$ and $M$ a finitely generated $R$-module. Let $S$ be a multiplicative
subset of $R$ with $S\cap \fa=\emptyset$.
\begin{enumerate}
\item[(i)] If $M$ is $\fa$-relative Cohen-Macaulay, then $S^{-1}M$ is $\fa S^{-1}R$-relative Cohen-Macaulay.
\item[(ii)] If $M$ is $\fa$-relative maximal Cohen-Macaulay and $S\cap (\Ann_RM+\fa)=\emptyset$, then $S^{-1}M$ is
$\fa S^{-1}R$-relative maximal Cohen-Macaulay.
\item[(iii)] If $M$ is $\fa$-relative Gorenstein, then $S^{-1}M$ is $\fa S^{-1}R$-relative Gorenstein.
\item[(iv)] If $M$ is $\fa$-relative regular, then $S^{-1}M$ is $\fa S^{-1}R$-relative regular.
\end{enumerate}
\end{corollary}

\begin{prf} (i), (iii), and (iv) are clear by Theorem \ref{4.2}.

(ii) If $S\cap (\Ann_RM+\fa)=\emptyset$, then $S^{-1}M\neq (S^{-1}\fa)S^{-1}M$, and so the claim follows by Theorem \ref{4.2}(ii).
\end{prf}

In Theorem \ref{4.2}, we may wonder whether the assumption that $f$ is flat can be replaced with the weaker assumption
that $f$ has finite flat dimension. Here is an example showing that this is not the case.

\begin{example}\label{4.4} Let $R=\Bbbk[x,y]$ be a polynomial ring over a field $\Bbbk$ and set $T=R/\langle xy,y^2
\rangle $. Let $f:R\to T$ be the natural ring epimorphism. As $\pd_RT=2$, we have $\fd_RT<\infty$. We set $\fa=\langle x \rangle$
and $M=R$. Then $M$ is $\fa$-relative regular with $\cd(\fa,R)=1$. Thus, $M$ is $\fa$-relative Gorenstein, $\fa$-relative maximal
Cohen-Macaulay, and $\fa$-relative Cohen-Macaulay.

On the other hand, we claim that $T\otimes_RM=T$ is not $\fa T$-relative Cohen-Macaulay. To this end, we show $\grade(\fa T,T)=0$ and
$\cd(\fa T,T)=1$. Observe that $\Ass_RT=\{\langle y \rangle, \langle x,y \rangle \}.$  Clearly, $\cd(\fa,R/\langle y\rangle )=1$ and $\cd(\fa,R/\langle x,y\rangle )=0$. Hence, $$\cd(\fa T,T)=\cd(\fa,T)=\max\{\cd(\fa,R/\langle y\rangle), \cd(\fa,R/\langle x,y
\rangle)\}=1.$$ As $\langle x,y\rangle \in \Ass_RT$ and it containing $\fa$, we have $$\grade(\fa T,T)=\grade(\fa,T)=0.$$
So, $T$ is not $\fa T$-relative Cohen-Macaulay. Thus, $T$ is not $\fa T$-relative maximal Cohen-Macaulay, not $\fa T$-relative
Gorenstein, and not $\fa T$-relative regular.
\end{example}

The following example indicates that in Theorem \ref{4.2}, the assumption that $f$ is flat can not be replaced with the assumption
that $f$ is pure.

\begin{example}\label{4.5} Let $\fa$ be a proper ideal of $R$ and $M$ a finitely generated $R$-module. Assume that $R$ is $\fa$-relative
regular and $M$ is not $\fa$-relative maximal Cohen-Macaulay. Let $T=R\ltimes M$ denote the trivial extension of $R$ by $M$. By Remark \ref{3.10}(ii), the canonical injection $\lambda:R\lo T$ is pure. Set $N=R$. Then $N$ is $\fa$-relative regular, and so it is $\fa$-relative Gorenstein, $\fa$-relative maximal Cohen-Macaulay and $\fa$-relative Cohen-Macaulay. On the other hand, $T\otimes_RN\cong T$ is not $\fa T$-relative Cohen-Macaulay by Proposition \ref{3.12}(i). Thus, $T\otimes_RN$ is not $\fa T$-relative maximal Cohen-Macaulay, not
$\fa T$-relative Gorenstein, and not $\fa T$-relative regular.
\end{example}

\begin{proposition}\label{4.6} Let $f:R\lo T$ be a faithfully flat ring homomorphism, $\fa$ a proper ideal of $R$, and $M$ a finitely
generated $R$-module.
\begin{enumerate}
\item[(i)] $M$ is $\fa$-relative Cohen-Macaulay if and only if $T\otimes_RM$ is $\fa T$-relative Cohen-Macaulay.
\item[(ii)] $M$ is $\fa$-relative maximal Cohen-Macaulay if and only if $T\otimes_RM$ is $\fa T$-relative maximal Cohen-Macaulay.
\item[(iii)] $M$ is $\fa$-relative Gorenstein if and only if $T\otimes_RM$ is $\fa T$-relative Gorenstein.
\end{enumerate}
\end{proposition}

\begin{prf} Since every faithfully flat ring homomorphism is both flat and pure, (i) is immediate by Theorems \ref{3.3}(i) and
\ref{4.2}(i).

(ii) As $f$ is faithfully flat, for every $R$-module $X$, it follows that $X=0$ if and only if $T\otimes_RX=0$. In particular,
$M\neq \fa M$ if and only if $T\otimes_RM\neq (\fa T)(T\otimes_RM)$. Now, the assertion follows by Theorems \ref{3.3}(ii) and
\ref{4.2}(ii).

(iii) Assume that $M=\fa M$. Then, $T\otimes_RM=(\fa T)(T\otimes_RM)$. These equalities yield that $\Ext_R^i(R/\fa,M)=0$ and 
$\Ext_T^i(T/\fa T,T\otimes_RM)=0$ for all $i\geq 0$, respectively. So, we may and do assume that $M\neq \fa M$. Then, as we 
observed in the proof of (ii), $T\otimes_RM\neq (\fa T)(T\otimes_RM)$. Thus, in view of Theorem \ref{2.4}(ii) and (ii), we may 
assume that $M$ is $\fa$-relative maximal Cohen-Macaulay and $T\otimes_RM$ is $\fa T$-relative maximal Cohen-Macaulay. Now, Lemma 
\ref{4.1}(ii) yields that $\cd(\fa,R)=\cd(\fa T,T)$.

As $f$ is flat, there is a natural $T$-isomorphism $$\Ext_T^i(T/\fa T,T\otimes_RM)\cong \Ext_R^i(R/\fa,M)\otimes_RT$$ for all
$i\geq 0$.  For each $i\geq 0$, in view of the faithful flatness of $f$, it turns out that $\Ext_R^i(R/\fa,M)=0$ if and only
if $\Ext_T^i(T/\fa T,T\otimes_RM)=0$. Therefore, $M$ is $\fa$-relative Gorenstein if and only if $T\otimes_RM$ is $\fa T$-relative
Gorenstein.
\end{prf}

Combining Theorems \ref{3.3}(iii) and \ref{4.2}(iv) yield the following:

\begin{corollary}\label{4.7} Let $f:R\lo T$ be a faithfully flat ring homomorphism, $\fa$ an ideal of $R$ contained in its
Jacobson radical and $M$ a finitely generated $R$-module. Assume that $\mu(\fa)=\cd(\fa,R)$. Then $M$ is $\fa$-relative regular
if and only if $T\otimes_RM$ is $\fa T$-relative regular.
\end{corollary}

Next, we record the following immediate three corollaries of Proposition \ref{4.6} and Theorem \ref{4.2}(iv).

\begin{corollary}\label{4.8} Let $\fa\subseteq \fb$ be two ideals of $R$ contained in its Jacobson radical and $M$ a finitely generated
$R$-module. Let $T$ denote the $\fb$-adic completion of $R$.
\begin{enumerate}
\item[(i)]  $M$ is $\fa$-relative Cohen-Macaulay if and only if $T\otimes_RM$ is $\fa T$-relative Cohen-Macaulay.
\item[(ii)]  $M$ is $\fa$-relative maximal Cohen-Macaulay if and only if $T\otimes_RM$ is $\fa T$-relative maximal Cohen-Macaulay.
\item[(iii)]  $M$ is $\fa$-relative Gorenstein if and only if $T\otimes_RM$ is $\fa T$-relative Gorenstein.
\item[(iv)] If $M$ is $\fa$-relative regular, then $T\otimes_RM$ is $\fa T$-relative regular.
\end{enumerate}
\end{corollary}

In particular, we have:

\begin{corollary}\label{4.9} Let $\fa$ be a proper ideal of a local ring $(R,\fm)$ and $M$ a finitely generated $R$-module.
\begin{enumerate}
\item[(i)] $M$ is $\fa$-relative Cohen-Macaulay if and only if $\widehat{M}$ is $\fa \widehat{R}$-relative Cohen-Macaulay.
\item[(ii)] $M$ is $\fa$-relative maximal Cohen-Macaulay if and only if $\widehat{M}$ is $\fa \widehat{R}$-relative maximal
Cohen-Macaulay.
\item[(iii)] $M$ is $\fa$-relative Gorenstein if and only if $\widehat{M}$ is $\fa \widehat{R}$-relative Gorenstein.
\item[(iv)] If $M$ is $\fa$-relative regular, then $\widehat{M}$ is $\fa \widehat{R}$-relative regular.
\end{enumerate}
\end{corollary}

\begin{corollary}\label{4.10} Let $\fa$ be a proper ideal of $R$ and $M$ a finitely generated $R$-module. Let $x_1,\dots, x_n$
be $n$ indeterminates over $R$.
\begin{enumerate}
\item[(i)]  $M$ is $\fa$-relative Cohen-Macaulay if and only if $M[x_1,\dots, x_n]$ is $\fa R[x_1,\dots, x_n]$-relative
Cohen-Macaulay if and only if $M[[x_1,\dots, x_n]]$ is $\fa R[[x_1,\dots, x_n]]$-relative
Cohen-Macaulay.
\item[(ii)]  $M$ is $\fa$-relative maximal Cohen-Macaulay if and only if $M[x_1,\dots, x_n]$ is
$\fa R[x_1,\dots, x_n]$-relative maximal Cohen-Macaulay if and only if $M[[x_1,\dots, x_n]]$ is
$\fa R[[x_1,\dots, x_n]]$-relative maximal Cohen-Macaulay.
\item[(iii)]  $M$ is $\fa$-relative Gorenstein if and only if $M[x_1,\dots, x_n]$ is $\fa R[x_1,\dots, x_n]$-relative
Gorenstein if and only if $M[[x_1,\dots, x_n]]$ is $\fa R[[x_1,\dots, x_n]]$-relative
Gorenstein.
\item[(iv)] If $M$ is $\fa$-relative regular, then $M[x_1,\dots, x_n]$ and $M[[x_1,\dots, x_n]]$ are $\fa R[x_1,\dots, x_n]$-relative
regular and $\fa R[[x_1,\dots, x_n]]$-relative regular, respectively.
\end{enumerate}
\end{corollary}

\begin{proposition}\label{4.11} Let $\fa$ be an ideal of $R$ and $M$ a finitely generated $R$-module with $M\neq \fa M$. Let $\underline{x}=x_1,x_2, \dots, x_n$ be a part of an $\fa$-s.o.p of $M$ which is an $M$-regular sequence. Set $\overline{R}=
R/\langle \underline{x}\rangle$ and $\overline{M}=M/\langle \underline{x} \rangle M$.  The following are equivalent:
\begin{enumerate}
\item[(i)] the $R$-module $M$ is $\fa$-relative Cohen-Macaulay;
\item[(ii)] the $R$-module $\overline{M}$ is $\fa$-relative Cohen-Macaulay;
\item[(iii)] the $\overline{R}$-module $\overline{M}$ is $\fa \overline{R}$-relative Cohen-Macaulay.
\end{enumerate}
\end{proposition}

\begin{prf} As $\underline{x}$ is an $M$-regular sequence, it follows that $$\grade(\fa \overline{R},\overline{M})=\grade(\fa,\overline{M})=\grade(\fa,M)-n.$$  On the other hand, as $\underline{x}$ is a part of an
$\fa$-s.o.p of $M$, Theorem \ref{2.2} yields that $$\cd(\fa \overline{R},\overline{M})=\cd(\fa,\overline{M})=\cd(\fa,M)-n.$$
The assertion is immediately implied by these equalities.
\end{prf}

\begin{proposition}\label{4.12} Let $\fa$ be a proper ideal of $R$ and $M$ a finitely generated $R$-module. Let $\underline{x}=x_1,
x_2, \dots, x_n$ be a part of an $\fa$-s.o.p of $R$, and set $\overline{R}=R/\langle \underline{x} \rangle$ and $\overline{M}=M/\langle
\underline{x} \rangle M$.
\begin{enumerate}
\item[(i)] Assume that $\underline{x}$ is an $M$-regular sequence. Then $M$ is $\fa$-relative maximal Cohen-Macaulay if and only
if the $\overline{R}$-module $\overline{M}$ is $\fa \overline{R}$-relative maximal Cohen-Macaulay.
\item[(ii)] Assume that $\underline{x}$ is both an $M$-regular sequence and an $R$-regular sequence. Then $M$ is $\fa$-relative
Gorenstein if and only if the $\overline{R}$-module $\overline{M}$ is $\fa \overline{R}$-relative Gorenstein.
\end{enumerate}
\end{proposition}

\begin{prf} (i) As $x_1, x_2, \dots, x_n$ is a part of an $\fa$-s.o.p of $R$, Theorem \ref{2.2} yields that $$\cd(\fa \overline{R},\overline{R})=\cd(\fa,\overline{R})=\cd(\fa,R)-n.$$ Also, as $\underline{x}$ is an $M$-regular sequence, it follows
that $$\grade(\fa \overline{R},\overline{M})=\grade(\fa,\overline{M})=\grade(\fa,M)-n.$$ Thus, $M$ is an $\fa$-relative maximal
Cohen-Macaulay if and only if $\overline{M}$ is $\fa \overline{R}$-relative maximal Cohen-Macaulay.

(ii) We have $$\cd(\fa \overline{R},\overline{R})=\cd(\fa,R)-n,$$ as we saw in the proof (i). As $\underline{x}$ is an $M$-regular
sequence, $\grade(\fa,M)\geq n$, and so $\Ext_R^i(R/\fa,M)=0$ for all $i<n$. Clearly, $R/\fa\cong \overline{R}/\fa \overline{R}$
and by the assumption $\underline{x}$ is both an $M$-regular sequence and an $R$-regular sequence. Thus, \cite[\S 18, Lemma 2(i)]{M}
implies the $R$-isomorphism $$\Ext_{\overline{R}}^i(\overline{R}/\fa \overline{R},\overline{M})\cong \Ext_R^{i+n}(R/\fa,M)$$ for all
$i\geq 0$. Hence, $\Ext_R^i(R/\fa,M)=0$ for all $i\neq \cd(\fa,R)$ if and only if $\Ext_{\overline{R}}^i(\overline{R}/\fa \overline{R},\overline{M})=0$ for all $i\neq \cd(\fa \overline{R},\overline{R})$. Thus, $M$ is $\fa$-relative Gorenstein if and only
if $\overline{M}$ is $\fa \overline{R}$-relative Gorenstein.
\end{prf}

\begin{corollary}\label{4.13} Let $\fa$ be a proper ideal of $R$ and $\underline{x}=x_1, x_2, \dots, x_n$ a part of an $\fa$-s.o.p
of $R$. Set $\overline{R}=R/\langle \underline{x} \rangle$.
\begin{enumerate}
\item[(i)] Assume that $\underline{x}$ is an $R$-regular sequence. Then $R$ is $\fa$-relative Gorenstein if and only if the ring
$\overline{R}$ is $\fa \overline{R}$-relative Gorenstein.
\item[(ii)] Assume that $\fa$ is contained in its Jacobson radical of $R$. If $R$ is $\fa$-relative Gorenstein, then the ring
$\overline{R}$ is $\fa \overline{R}$-relative Gorenstein.
\end{enumerate}
\end{corollary}

\begin{prf} (i) is obvious by Proposition \ref{4.12}(ii).

(ii) Suppose that $R$ is $\fa$-relative Gorenstein. Then, $R$ is $\fa$-relative Cohen-Macaulay by Theorem \ref{2.4}. As $\underline{x}$
is a part of an $\fa$-s.o.p of $R$ and $\fa$ is contained in the Jacobson radical of $R$, Theorem \ref{2.5} asserts that $\underline{x}$
is an $R$-regular sequence. Now, (i) implies that the ring $\overline{R}$ is $\fa \overline{R}$-relative Gorenstein.
\end{prf}

\begin{proposition}\label{4.14} Let $\fa$ be an ideal of $R$, $M$ a finitely generated $R$-module with $M\neq \fa M$, and set
$c=\cd(\fa,R)$. Assume that $M$ is $\fa$-relative regular and $x_1, x_2,\dots, x_c$ is a generating set of $\fa$. Let $1\leq
n\leq c$ be an integer, and set $\overline{R}=R/\langle x_1, x_2,\dots, x_n\rangle$ and $\overline{M}=M/\langle x_1, x_2,\dots,
x_n \rangle M.$ Then the $\overline{R}$-module $\overline{M}$ is $\fa \overline{R}$-relative regular, provided that either
\begin{enumerate}
\item[(i)] $x_1, x_2,\dots, x_c$ form both an $M$-regular sequence and an $R$-regular sequence; or
\item[(ii)] $\fa$ is contained in the Jacobson radical of $R$.
\end{enumerate}
\end{proposition}

\begin{prf} (i) Clearly, $\fa \overline{R}=\langle \overline{x}_{n+1}, \overline{x}_{n+2}, \dots, \overline{x}_{c} \rangle$ and $\overline{x}_{n+1}, \overline{x}_{n+2}, \dots, \overline{x}_{c}$ is both an $\overline{M}$-regular sequence and an
$\overline{R}$-regular sequence. As $\fa=\langle x_1, x_2,\dots, x_c\rangle$, it follows that $x_1, x_2,\dots, x_c$ is an
$\fa$-s.o.p of $R$. Now, Theorem \ref{2.2} implies that $\cd(\fa \overline{R},\overline{R})=c-n$, and so $\overline{M}$ is
$\fa \overline{R}$-relative regular.

(ii) Since $M$ is $\fa$-relative regular and $M\neq \fa M$,  the ring $R$ is also $\fa$-relative regular and $\cd(\fa,M)=\cd(\fa,R)$.
By Theorem \ref{2.4}, both $M$ and $R$ are $\fa$-relative Cohen-Macaulay. Since $\fa=\langle x_1, x_2,\dots, x_c \rangle,$ it follows
that $x_1, x_2, \dots, x_c$ is both an $\fa$-s.o.p of $M$ and an $\fa$-s.o.p of $R$. Now, Theorem \ref{2.5} yields that $x_1, x_2,
\dots, x_c$ is both an $M$-regular sequence and an $R$-regular sequence, and so the claim follows by (i).
\end{prf}

\begin{corollary}\label{4.15} Let $\fa$ be a proper ideal of $R$, and set $c=\cd(\fa,R)$. Assume that $R$ is $\fa$-relative regular
and $x_1, x_2,\dots, x_c$ is a generating set of $\fa$. Let $1\leq n\leq c$ be an integer and set $\overline{R}=R/\langle x_1, x_2,
\dots, x_n\rangle$. Then $\overline{R}$ is $\fa \overline{R}$-relative regular, provided that either
\begin{enumerate}
\item[(i)] $x_1, x_2,\dots, x_c$ form an $R$-regular sequence; or
\item[(ii)] $\fa$ is contained in the Jacobson radical of $R$.
\end{enumerate}
\end{corollary}

%%%%%%%%%%%%%%%%%%%%%%%%%%%%%%%%%%%%%%%%%%%%%%%%%%%

\end{document}